\documentclass[12pt]{article}
\usepackage{amssymb}
\def\qed{\hfill$\Box$\par}
\def\a{\alpha}

\def\d{\delta}

\def\g{\gamma}

\def\l{\lambda}

\def\si{\sigma}

\def\Vir{\hbox{Vir}}
\def\sc{\scriptstyle}
\def\ssc{\scriptscriptstyle}
\def\dis{\displaystyle}
\def\cl{\centerline}
\def\ll{\leftline}

\def\nl{\newline}

\def\ul{\underline}
\def\wt{\widetilde}

\def\Lra{\Leftrightarrow}

\def\hs{\hspace*}
\def\vs{\vspace*}

\def\ni{\noindent}
\def\ptl{\partial}
\def\ptla{\mbox{$\ptl\over\ptl\l$}}
\def\hi{\hangindent}
\def\ha{\hangafter}

\def\Z{\mathbb{Z}}

\def\C{\mathbb{C}}

\def\NN{{\cal N}_q}\def\ll{{\ul\l}}\def\SS{{\cal S}_q}\def\nn{{\ul
n}}\def\mm{{\ul m}}\def\kk{{\ul k}}%
\begin{document}
%
%
\ni{\large\bf Low dimensional cohomology 
 of general conformal algebras
$gc_N$} 
\cl{(appeared in {\em J. Math. Phys.} {\bf45} (2004) 509--524)}
\\
\begin{minipage}[rm]{6truein}
\parskip .05 truein \baselineskip 4pt \lineskip 4pt
\ni Yucai Su
\par%
{\small\it Department of Mathematics,
                Shanghai Jiaotong University, Shanghai 200030,
                P.~R.~China
                 {\rm and}
Department of Mathematics,
                Harvard University, Cambridge, MA 02138, USA
\nl
                Email: ycsu@sjtu.edu.cn
}
\par\ \vs{-5pt}\par\ni
{\small 
We compute the low dimensional cohomologies $\wt H^q(gc_N,\C),$
$H^q(gc_N,\C)$
 of the
infinite rank general Lie conformal algebras $gc_N$ with trivial
coefficients for $q\le3,\,N=1$ or $q\le2,\,N\ge2$. We also prove
that the cohomology of $gc_N$ with coefficients in its natural
module is trivial, i.e., $H^*(gc_N,\C[\ptl]^N)=0$; thus partially
solve an open problem of Bakalov-Kac-Voronov in [{\it
Comm.~Math.~Phys.,} {\bf200} (1999), 561-598].}
\end{minipage}
\par\
\vs{-5pt}\par \ni {\bf
             I. \ INTRODUCTION
} \par
The notion of a conformal algebra, introduced by Kac in Ref.~12,
encodes an axiomatic description of the operator product expansion
of chiral fields in conformal field theory. Conformal algebras
play important roles in quantum field theory and vertex operator
algebras (e.g.~Ref.~12), whose study has drawn much attention in
the literature (e.g.~Refs.~1--6, 10, 12--14 and 20--23). As is
pointed out in Ref.~2, on one hand, it is an adequate tool for the
study of infinite-dimensional Lie algebras satisfying the locality
property (cf.~Refs.~5, 12 and 14). On the other hand, conformal
modules over a conformal algebra $R$ correspond to conformal
modules over the associated Lie algebra Lie$\,R$ (cf.~Ref.~3). The
main examples of Lie algebras Lie$\,R$ are the Lie algebras
``based'' on the punctured complex plane $\C^\times$, namely the
Lie algebra {\sl Vect}$\,\C^\times$ of vector fields on
$\C^\times$ (the Virasoro algebra) and the Lie algebra of maps of
$\C^\times$ to a finite-dimensional Lie algebra (the loop
algebra). Their irreducible conformal modules are the spaces of
densities on $\C^\times$ and loop modules, respectively
(cf.~Ref.~3). Since complete reducibility does not hold in this
case (cf.~Refs.~4 and 9), one may expect that their cohomology
theory is very interesting and important (cf.~Ref.~2), just as the
cohomology theory of Lie algebras has played important roles in
the structure and representation theories of Lie algebras
(cf.~Refs.~7--9, 11 and 15--19).
\par
A general theory of cohomology of Lie conformal algebras was
established by Bakalov, Kac and Voronov in Ref.~2. They also
computed the cohomologies for the finite simple Lie conformal
algebras. However the problem for the general Lie conformal
algebra $gc_N$, which is an infinite Lie conformal algebra,
remains open. It is well-known that the general Lie conformal
algebra $gc_N$ plays the same important role in the theory of Lie
conformal algebras as the general Lie algebra $gl_N$ does in the
theory of Lie algebras: any module $M = \C[\ptl]^N$ over a Lie
conformal algebra $R$ is obtained via a homomorphism $R\to gc_N$
(cf.~Refs.~5 and 12), thus the study of Lie conformal algebras
$gc_N$ has drawn some authors' attentions (cf.~Refs.~1, 2, 6, 13
and 14). It seems to us that the computation of cohomology of
$gc_N$
is important.%
\vs{-3pt}\par%
In this paper, we compute the low dimensional basic cohomologies
$\wt H^q(gc_N,\C)$ and \vs{-3pt}the reduced cohomology $\wt
H^q(gc_N,\C)$ of $gc_N$ with trivial coefficients for
$q\!\le\!3,\,N\!=\!1$ or $q\!\le\!2,\,N\!\ge\!2$. We also prove
that the cohomology of $gc_N$ with coefficients in its natural
module is trivial, i.e., $H^*(gc_N,\C[\ptl]^N)=0$; thus partially
solve an open problem 
 in Ref.~2. \par%
In Section II, we shall recall definitions of conformal algebras,
their modules and cohomology, and present the main theorem of this
paper (Theorem 2.5). Sections III and IV are devoted to the proof
of the main theorem.
\par\ \vs{-5pt}\par\ni{\bf II. \ NOTATIONS AND MAIN RESULTS}
\par%
We shall briefly recall definitions of conformal algebras, their
modules and cohomology. More details can be
found in, say, Ref.~2. \par%
{\it Definition 2.1:} A {\it Lie conformal algebra} is a
$\C[\ptl]$-module $A$ with a $\l$-bracket $[a{\ssc\,}_\l b]$ which
defines a linear map $A\times A\to A[\l],$ where $A[\l] =
\C[\l]\otimes A$ is the space of polynomials of $\l$ with
coefficients in $A$,
\vs{-6pt}satisfying: %
$$%
\begin{array}{ll}[\ptl a{\ssc\,}_\l b] =-\l[a{\ssc\,}_\l b],
\ \ [a{\ssc\,}_\l\ptl b]=(\ptl+\l)[a{\ssc\,}_\l b]
&\mbox{ \ (conformal sesquilinearity),}\vs{4pt}\\ %
{}[a{\ssc\,}_\l b]=-[b{\ssc\,}_{-\l-\ptl}a]&\mbox{ \ (skew-symmetry)},\vs{4pt}\\ %
{}[a{\ssc\,}_\l[b{\ssc\,}_\mu c]]=[[a{\ssc\,}_\l
b]_{\l+\mu}c]+[b{\ssc\,}_\mu[a{\ssc\,}_\l c]]&\mbox{ \ (Jacobi
identity)},
\end{array}%
\vs{-6pt}\eqno\begin{array}{r}\vs{4pt}(2.1)\!\!\\ \vs{4pt}(2.2)\!\!\\
(2.3)\!\!\end{array}
$$%
for $a,b,c\!\in\! A$. A subset $S\!\subset\! A$ is called a {\it
generating set} if $S$ generates $A$ as a $\C[\ptl]$-module. If
there exists a finite generating set, then $A$ is called {\it
finite}. Otherwise, it is called {\it infinite}.
\qed%
There is a similar notion of associative conformal algebras, which
we shall not introduce in this paper. Below we shall only work
with Lie conformal algebras, thus we shorten the term ``Lie
conformal algebra'' to ``conformal algebra''. The simplest
nontrivial conformal algebra is the {\it Virasoro conformal
algebra} \Vir, which is a rank one free $\C[\ptl]$-module
generated by a symbol $L$ such \vs{-4pt}that
$$%
\Vir=\C[\ptl]L,\;\;\;[L{\ssc\,}_\l L]=(\ptl +2\l)L.%
\vs{-4pt}\eqno(2.4)$$%
Note that using (2.1), it suffices to define $\l$-brackets on a
generating set. Let $N\ge1$ be an integer. The {\it general
conformal algebra $gc_N$} can be defined (see, e.g., Ref.~14) as
an infinite rank
free $\C[\ptl]$-module with a generating \vs{-4pt}set%
$$S_N=\{J^n_A\,|\,n\in\Z_+,\,A\in gl_N\},%
\vs{-4pt}\eqno(2.5)$$%
 where $gl_N$ is the space of $N\times N$ matrices (note that the set $S_N$ is not
 $\C$-linearly independent, for example, $J^m_{aA}=aJ^m_A$ for $a\in\C$), such that the \vs{-7pt}$\l$-bracket is defined by
$$%
[J^m_A{\ssc\,}{\ssc\,}_\l
J^n_B]=\sum_{s=0}^m\biggl(\!\begin{array}{c}m\\
s\end{array}\!\biggr)(\l+\ptl)^sJ^{m+n-s}_{AB}-\sum_{s=0}^n\biggl(\!\begin{array}{c}n\\
s\end{array}\!\biggr)(-\l)^s
J^{m+n-s}_{BA},%
\vs{-7pt}\eqno(2.6)$$%
for $m,n\in\Z_+,\,A,B\in gl_N$, where
$(^m_{{\ssc\,}s})={m(m-1)\cdots(m-s+1)\over s!}$ if $s\ge0$ and
$(^m_{{\ssc\,}s})=0$ otherwise, is the binomial coefficient.
\par%
{\it Definition 2.2:} A {\it module} over a conformal algebra $A$
is a $\C[\ptl]$-module $M$ with a $\l$-action $a{\ssc\,}_\l v$
which defines a map $A\times M\to M[[\l]]$, where $M[[\l]]$ is the
set of formal power series of $\l$ with coefficients in $M$, such
that
$$%
\begin{array}{l}
a{\ssc\,}_\l(b{\ssc\,}_\mu v)-b{\ssc\,}_\mu(a{\ssc\,}_\l v)=[a{\ssc\,}_\l b]_{\l+\mu}v,\vs{4pt}\\ %
(\ptl a){\ssc\,}_\l v=-\l a{\ssc\,}_\l v,\;\;\;a{\ssc\,}_\l(\ptl v)=(\ptl+\l)a{\ssc\,}_\l v,%
\end{array}%
\eqno\begin{array}{r}\vs{4pt}(2.7)\!\!\\ (2.8)\!\!\end{array}
$$%
for $a,b\in A,\,v\in M$. If $a{\ssc\,}_\l v\in M[\l]$ for all
$a\in A,\,v\in M$, then the $A$-module $M$ is called {\it
conformal}. If $M$ is finitely generated over $\C[\ptl]$, then $M$
is simply called {\it finite}.
\par%
Below we shall only consider ``conformal modules'', thus we drop
the word ``conformal'' and simply call a ``conformal module'' a
``module''. Clearly, the one-dimensional vector space $\C$ can be
defined as a module (called a {\it trivial module}) over any
conformal algebra $A$ with both the action of $\ptl$ and the
action of $A$ being zero. Furthermore, for $a\in\C,\,a\ne0$, one
can define a $\C[\ptl]$-module $\C_a$, which is the
one-dimensional vector space $\C$ such that $\ptl v=av$ for
$v\in\C_a$. Then $\C_a$ becomes an $A$-module with trivial action
of $A$.
\par
Let $\a\in\C$. The space $\C^N[\ptl]$ (a rank $N$ free
$\C[\ptl]$-module) can be defined as a $gc_N{\ssc\,}$-module with\vs{-4pt} $\l$-action%
$$%
J^m_A{\ssc\,}_\l v=(\ptl+\l+\a)^m Av\mbox{ \ \ for \ }A\in
gl_N,\,m\in\Z_+,\,v\in\C^N,%
\vs{-4pt}\eqno(2.9)
$$%
(cf.~the statement after (2.4)). We denote this module by
$\C_\a^N[\ptl]$. When $\a=0$, the module $\C^N[\ptl]=\C_0^N[\ptl]$
is called the {\it natural module of $gc_N$.}
\par%
{\it Definition 2.3:} Let $q\in\Z_+$. A {\it $q$-cochain of a
conformal algebra A with coefficients in a module $M$} is a
$\C$-linear map\vs{-4pt} $\g:A^{\otimes q}\to M[\l_1,...,\l_q]$,%
$$%
\g(a_1\otimes\cdots\otimes a_q)=\g_{\l_1,...,\l_q}(a_1,...,a_q),
\vs{-4pt}\eqno(2.10)$$%
satisfying%
\par\ni\hspace*{2ex}$%
\!\!\!\!\begin{array}{l}%
\g_{\l_1,...,\l_q}(a_1,...,\ptl a_i,...a_q)=-\l_i
\g_{\l_1,...,\l_q}(a_1,...,a_i,...,a_q)\mbox{ \ (conformal antilinearity)},\vs{4pt}\\%
\g_{\l_1,...,\l_{i-1},\l_{i+1},\l_i,\l_{i+2},...,\l_q}(a_1,...,a_{i-1},a_{i+1},a_i,a_{i+2},...,a_q)
\vs{2pt}\\
\ \ \ \ \ =-
\g_{\l_1,...,\l_{i-1},\l_i,\l_{i+1},\l_{i+2},...,\l_q}(a_1,...,a_{i-1},a_i,a_{i+1},a_{i+2},...,a_q)
\mbox{ \ (skew-symmetry)},
\end{array}$\hfill%
$\begin{array}{r}\vs{4pt}(2.11)\!\!\\
\vs{2pt}\ \\ (2.12)\!\!\end{array}
$\par\ni%
for $a_1,...,a_q\in A$ and all possible $i$. We let $A^{\otimes
0}=\C$, so
that a $0$-cochain $\g$ is simply an element of $M$. %
\qed%
We define a {\it differential $\,d\,$ of a cochain $\g$} as
follows:
$$%
\begin{array}{ll}%
&\dis(d\g)_{\l_1,...,\l_{q+1}}(a_1,...,a_{q+1})\vs{2pt}\\
=\!\!\!\!&\dis
 \sum_{i=1}^{q+1}(-1)^{i+1}a_i{\ssc\,}_{\l_i}
\g_{\l_1,...,\hat\l_i,...,\l_{q+1}}(a_1,...,\hat
a_i,...,a_{q+1})\vs{0pt}\\ %
&+\dis\sum_{1\le i<j\le
q+1}(-1)^{i+j}\g_{\l_i+\l_j,\l_1,...,\hat\l_i,...,\hat\l_j,...,\l_{q+1}}
([a_i{\ssc\,}_{\l_i}a_j],a_1,...,\hat a_i,...,\hat
a_j,...,a_{q+1}),
\end{array}%
\vs{-7pt}\eqno(2.13)$$%
where, $\g$ is extended linearly over the polynomials in $\l_i$,
and where, the symbol $\hat{\ }$ means the element below it is
missing. In particular,
$$(d\g){\ssc\,}_\l(a)=a{\ssc\,}_\l\g\mbox{ \ if \ $\g\in M$ \ is a $0$-cochain}.%
\eqno(2.14)$$%
By Ref.~2, the operator $d$ preserves the space of cochains and
$d^2=0$, so that the cochains form a complex, which will be
denoted by $\wt C^*=\wt C^*(A,M)=\oplus_{q\in\Z_+}\wt C^q(A,M)$,
and called the {\it basic complex}
for the $A$-module $M$.%
\par%
Define the structure of a $\C[\ptl]$-module on $\wt C^*(A,M)$
\vs{-9pt}by
$$%
(\ptl\g)_{\l_1,...,\l_q}(a_1,...,a_q)=
(\ptl_M+\sum_{i=1}^q\l_i)\g_{\l_1,...,\l_q}(a_1,...,a_q),%
\vs{-9pt}\eqno(2.15)
$$%
where $\ptl_M$ denotes the action of $\ptl$ on $M$. Then
$d\ptl=\ptl d$ (see Ref.~2) and so the graded subspace $\ptl\wt
C^*\subset\wt C^*$ forms a subcomplex. Define the quotient complex
$C^*=C^*(A,M)=\wt C^*(A,M)/\ptl\wt
C^*(A,M)=\oplus_{q\in\Z_+}C^q(A,M)$, called the {\it reduced
complex}.%
\par%
{\it Definition 2.4:} The {\it basic cohomology $\wt H^*(A,M)$ of
a conformal algebra $A$ with coefficients in a module $M$} is the
cohomology of the basic complex $\wt C^*$. The {\it (reduced)
cohomology $H^*(A,M)$} is the cohomology of the reduced complex
$C^*$.%
\qed%
Note that the basic cohomology $\wt H^*(A,M)$ is naturally a
$\C[\ptl]$-module, whereas the reduced cohomology $H^*(A,M)$ is a complex vector space.%
\par%
The main results of this paper is the following theorem.
\par%
{\bf Theorem 2.5:} {\it (1) For the general conformal algebra
$gc_1$, we \vs{-4pt}have%
$${\rm dim\,}\wt H^q(gc_1,\C)=\left\{\begin{array}{ll}1&\mbox{if
\ $q=0$ or $3$,}\\0&\mbox{if \ $q=1$ or $2$,}\end{array}\right.
\vs{-4pt}\eqno(2.16)$$%
\vs{-4pt}and%
$${\rm dim\,}
H^q(gc_1,\C)=\left\{\begin{array}{ll}1&\mbox{if \  $q=0,\,2$ or
$3$,}\\0&\mbox{if \ $q=1$;}\end{array}\right.%
\eqno(2.17)$$%
\vs{-4pt}\par%
(2) Equations (2.16) and (2.17) also hold for the general
conformal algebra $gc_N$ if $q\le2$;
%
\vskip -2pt\par%
(3) $H^*(gc_N,\C_a)=0$ if $a\ne0$;
\vs{-2pt}\par%
(4) $H^*(gc_N,\C_\a^N[\ptl])=0$ for $\a\in\C$. Furthermore, for
any $gc_N{\ssc\,}$-module $M$ which is freely generated over
$\C[\ptl]$ such that there exists nonzero $c\in\C$ satisfying
$J_I^0{\ssc\,}_\l v|_{\l=0}=cv$ for $v\in M$, where $I$ is the
$N\times N$ identity matrix, we have $H^*(gc_N,M)=0$.
}\vs{3pt}\par%
{\it Remark 2.6:} (1) Equations (2.16) and (2.17) show that the
cohomologies $\wt H^q(gc_1,\C),$ $H^*(gc_1,\C)$, $q\!\le\!3$, of
the general conformal algebra $gc_1$ with trivial coefficients are
isomorphic to those of the Virasoro conformal algebra with trivial
coefficients.
\vs{-2pt}\par%
(2) Theorem 2.5(2) in particular shows that there is a unique
nontrivial universal central extension of the general conformal
algebra $gc_N$, which agrees with that of the Lie algebra ${\cal
D}^N$ of $N\times N$ matrix differential operators on the circle
(cf.~Refs.~16 and 18. It is well-known that ${\cal D}^N$ is the
distribution Lie algebra associated with $gc_N$, cf.~Ref.~14). A
nontrivial reduced $2$-cocycle $\psi'$ of $gc_N$ is given in
(3.36), and the universal central extension $\wt{gc}_N$ of $gc_N$
corresponding to $\psi'$ is given \vs{-9pt}by
$$%
\begin{array}{ll}
[J^m_A{\ssc\,}_\l
J^n_B]=\!\!\!\!&\dis\sum_{s=0}^m\biggl(\!\begin{array}{c}m\\
s\end{array}\!\biggr)(\l+\ptl)^sJ^{m+n-s}_{AB}-\sum_{s=0}^n\biggl(\!\begin{array}{c}n\\
s\end{array}\!\biggr)(-\l)^s J^{m+n-s}_{BA}\\
&\dis+ (-1)^n{m!n!\over(m+n+1)!}{\rm tr}(AB)\l^{m+n+1}C,
\end{array}
\vs{-9pt}\eqno(2.18)$$%
where $C$ is a nonzero {\it central element} of $\wt{gc}_N$ (i.e.,
$[C{\ssc\,}_\l a]=[a{\ssc\,}_\l C]=0$ for all $a\in\wt{gc}_N$)
such that $\C C$ is a trivial $\C[\ptl]$-module. \vskip-2pt\par%
(3) In Theorem 2.5(4), note that if we define the $0$-bracket by
$[a{\ssc\,}_0b]=[a{\ssc\,}_\l b]|_{\l=0}$ for $a,b\in gc_N$, and
define the $0$-action of $gc_N$ on a module $M$ by
$a{\ssc\,}_0v=a{\ssc\,}_\l v|_{\l=0}$ for $a\in gc_N,\,v\in M$,
then $J^0_I$ is {\it central under $0$-bracket}, i.e.,
$[J_I^0{\ssc\,}_0a]=[a{\ssc\,}_0J_I^0]$ for $a\in gc_N$, and so
the $0$-action of $J^0_I$ on any indecomposable
$gc_N{\ssc\,}$-module $M$ is a scalar.
\qed%
We shall give the proof of Theorem 2.5 in the next two sections.
\par\vskip-5pt\ \par\ni{\bf III. \ PROOF OF THEOREM 2.5(2)-(4)}\par
We shall keep notations of the previous section. For a $q$-cochain
$\g\in\wt C^q(A,M)$, we call $\g$ a {\it $q$-cocycle} if $d\g=0$;
a {\it $q$-coboundary} or a {\it trivial $q$-cocycle} if there is
a \mbox{$(q-1)$-cochain} $\phi\in\wt C^{q-1}(A,M)$ such that
$\g=d\phi$. Two cochains $\g$ and $\psi$ are {\it equivalent} if
$\g-\psi$ is a coboundary. Denote by $\wt D^q(A,M)$ and by $\wt
B^q(A,M)$ the spaces of $q$-cocycles and $q$-coboundaries
respectively. Then by Definition
2.4, we have%
$$%
\wt H^q(A,M)=\wt D^q(A,M)/\wt B^q(A,M)=\{\mbox{equivalent classes
of
$q$-cocycles}\}.%
\eqno(3.1)
$$%
\par%
We shall divide the proof of Theorem 2.5(2)-(4) into several
lemmas (although we are unable to give the general result for
$gc_N$ in this paper, Lemmas 3.1-4 below may be helpful in
determining $\wt H^*(gc_N,\C)$ and $H^*(gc_N,\C)$ in the future).
\par%
First suppose $\g\in\wt C^q(gc_N,\C)$. Clearly, by (2.11), $\g$ is
uniquely determined by the right-hand side of (2.10) for
$a_1,...,a_q\in S_N$, where $S_N$ is defined in (2.5). We can
regard the right-hand side of (2.10) as a polynomial in
$\l_1,...,\l_q$.  For any fixed $p\in\Z$, we define a $\C$-linear
map $\g^{(p)}:gc_{{\sc N}}^{{\sc\otimes q}}\to\C[\l_1,...,\l_q]$
such that (2.11) holds for $\g^{(p)}$ and such that
$$%
\g^{(p)}(J_{A_1}^{n_1}\otimes\cdots\otimes
J_{A_q}^{n_q})=\g^{(p)}_{\l_1,...,\l_q}(J_{A_1}^{n_1},...,J_{A_q}^{n_q}),
\eqno(3.2)%
$$%
is a homogenous polynomial in $\l_1,...,\l_q$ consisting of all
monomials of total degree $p'$
 which
appear in $\g_{\l_1,...,\l_q}(J_{A_1}^{n_1},...,J_{A_q}^{n_q})$,
\vs{-9pt}where%
$$%
p'=p+\sum_{i=1}^qn_i.%
\vs{-9pt}\eqno(3.3)$$%
Then it is straightforward to see that $\g^{(p)}\in\wt
C^q(gc_N,\C)$ \vs{-7pt}and %
$$\g=\sum_{p\in\Z}\g^{(p)}.\vs{-9pt}\eqno(3.4)$$
Note that (3.4) is possibly an infinite sum, however for given
$J_{A_1}^{n_1},...,J_{A_q}^{n_q}\in S_N$, there are only finite
many $p$'s such that (3.2) is not zero; we call such a sum {\it
summable}. {}From (2.6), (2.11) and (2.13) (note that in (2.6), if
we informally regard the right-hand side as a polynomial in
$\l,\ptl,J_{AB},J_{BA}$, then it is a homogenous polynomial of the
total degree $m+n$; also note that (2.13) now takes the form such
that the first sum in the right-hand side is missing since $\C$ is
a trivial module and note from (2.11) that when we substitute
(2.6) into (2.13), $\ptl$ appeared in (2.6) can be replaced by
$-\l_i$ for some $i$), we immediately
obtain the following lemma. %
\par%
{\it Lemma 3.1:} {\it A $q$-cochain $\g\in\wt C^q(gc_N,\C)$ is a
$q$-cocycle (resp., $q$-coboundary) $\Lra$ all $\g^{(p)}$ are
$q$-cocycles (resp., $q$-coboundaries).}\hfill$\Box$.
\par
A $q$-cochain of the form $\g^{(p)}$ is called a {\it homogenous
$q$-cochain of degree $p$}.
\par%
Following Ref.~2, we define an operator $\tau_1:\wt C^q(gc_N,\C)\to \wt
C^{q-1}(gc_N,\C)$ as follows: If $q=0$, we set $\tau_1\g=0$;
otherwise, we set
$$%
(\tau_1\g)_{\l_1,...,\l_{q-1}}(a_1,...,a_{q-1})=
(-1)^{q-1}\ptla\g_{\l_1,...,\l_{q-1},\l}(a_1,...,a_{q-1},J)|_{\l=0},%
\eqno(3.5)$$%
for $a_1,...,a_{q-1}\in S_N$, where $J=J_I^1$ and $I$ is the
$N\times N$ identity matrix. Noting that by \vs{-7pt}(2.6),
$$%
[J^{n_i}_{A_i}{\ssc\,}_{\l_i}
J]=\sum_{s=1}^{n_i}\biggl(\!\begin{array}{c}n_i\\
s\end{array}\!\biggr)(\l_i+\ptl)^sJ^{n_i+1-s}_{A_i}-(-\l_i)J^{n_i}_{A_i},%
\vs{-7pt}\eqno(3.6)$$%
we obtain
$$%
\begin{array}{ll}
&((d\tau_1+\tau_1
d)\g^{(p)})_{\l_1,...,\l_q}(J_{A_1}^{n_1},...,J_{A_q}^{n_q})\vs{2pt}\\
=\!\!\!\!&\dis(-1)^q\ptla\sum_{i=1}^q(-1)^{i+q+1}\g^{(p)}_{\l_i+\l,\l_1,...,\hat\l_i,...,\l_q}
([J_{A_i}^{n_i}{\ssc\,}_{\l_i}J],J_{A_1}^{n_1},...,\hat
J_{A_i}^{n_i},...,J_{A_q}^{n_q})|_{\l=0}\vs{0pt}\\
=\!\!\!\!&\dis\ptla\sum_{i=1}^q\g^{(p)}_{\l_1,...,\l_{i-1},\l_i+\l,\l_{i+1},...,\l_q}
(J_{A_1}^{n_1},...,J_{A_{i-1}}^{n_{i-1}},[J_{A_i}^{n_i}{\ssc\,}_{\l_i}J],
J_{A_{i+1}}^{n_{i+1}},...,J_{A_q}^{n_q})|_{\l=0},
\end{array}
\vs{-4pt}\eqno(3.7)$$%
where the first equality follows from the fact that all terms
appearing in $d\tau_1\g^{(p)}$ are cancelled with the
corresponding terms in $\tau_1 d\g^{(p)}$ and the terms left are
all appearing in $\tau_1 d\g^{(p)}$ (cf.~(2.13)), the second
equality follows from (2.12). Note that for a polynomial $P$,
${\ptl P\over\ptl\l}|_{\l=0}$ is simply the coefficient of $\l^1$
in $P$. Now we substitute (3.6) into (3.7). By (2.11),
$(\l_i+\ptl)^s$ can be replaced by $(-\l)^s$. Since we only need
coefficients of $\l^1$, the terms with $s\ge2$ in (3.6) do not
contribute to the calculation. Thus
$[J_{A_i}^{n_i}{\ssc\,}_{\l_i}J]$ in (3.7) can be
replaced by $(\l_i-n_i\l)J_{A_i}^{n_i}$. Thus (3.7) is equal \vs{-9pt}to%
$$
\ptla\sum_{i=1}^q(\l_i-n_i\l)\g^{(p)}_{\l_1,...,\l_{i-1},\l_i+\l,\l_{i+1},...\l_q}
(J_{A_1}^{n_1},...,J_{A_q}^{n_q})|_{\l=0}
=p\g^{(p)}_{\l_1,...\l_q}(J_{A_1}^{n_1},...,J_{A_q}^{n_q}),%
\vs{-9pt}\eqno(3.8)%
$$%
which follows from (3.3) and the fact that for a homogenous
polynomial $P(\l_1,...,\l_q)$ of total degree $p'$, we
\vs{-9pt}have
$$\ptla\sum_{i=1}^q(\l_i-n_i\l)P(\l_1,...,\l_{i-1},\l_i+\l,\l_{i+1},...,\l_q)|_{\l=0}
=(p'-\sum_{i=1}^qn_i)P.%
\vs{-9pt}\eqno(3.9)$$%
{}From (3.7) and (3.8), we \vs{-5pt}obtain%
$$%
(d\tau_1+\tau_1 d)\g^{(p)}=p\g^{(p)}.\vs{-5pt}\eqno(3.10)
$$%
So if $d\g=0$, then (3.10) shows that
$\g'=\sum_{p\ne0}\g^{(p)}=d(\sum_{p\ne0}p^{-1}\tau_1\g^{(p)})$
(note that this is summable, cf.~the statement after (3.4)$\sc\,$)
is a
coboundary, and $\g-\g'=\g^{(0)}$. Thus, we obtain the following lemma.%
\par%
{\it Lemma 3.2:} {\it A $q$-cocycle in $\wt D^q(gc_N,\C)$ is
equivalent to a homogenous
$q$-cocycle of degree zero.}%
\qed%
Now suppose $\g$ is a homogenous $q$-cocycle of degree zero. For
$1\le j,k\le N$, denote by $E_{j,k}$ the $N\times N$ matrix with
entry $1$ at $(j,k)$ and $0$ otherwise. \vs{-4pt}Then %
$$%
S'_N=\{J^n_{E_{j,k}}\,|\,n\in\Z_+,\,1\le j,k\le N\},
\vs{-4pt}\eqno(3.11)$$%
is a free generating set of $gc_N$ over $\C[\ptl]$. Let
$h=\sum_{j=1}^N jE_{j,j}$. We define another operator $\tau_2:\wt
C^q(gc_N,\C)\to \wt C^{q-1}(gc_N,\C)$ as follows: We set
$\tau_2\g=0$ if $q=0$, otherwise we \vs{-3pt}set
$$%
(\tau_2\g)_{\l_1,...,\l_{q-1}}(a_1,...,a_{q-1})=
(-1)^{q-1}\g_{\l_1,...,\l_{q-1},0}(a_1,...,a_{q-1},J^0_h),%
\vs{-3pt}\eqno(3.12)$$%
for $a_1,...,a_{q-1}\in S_N$. Now note that by \vs{-9pt}(2.6),
$$%
[J^{n_i}_{E_{j_i,k_i}}{\ssc\,}_{\l_i}
J^0_h]=\sum_{s=0}^{n_i}\biggl(\!\begin{array}{c}n_i\\
s\end{array}\!\biggr)k_i(\l_i+\ptl)^sJ^{n_i-s}_{E_{j_i,k_i}}-j_iJ^{n_i}_{E_{j_i,k_i}}.%
\vs{-9pt}\eqno(3.13)$$%
Thus as discussion in (3.7) and (3.8), the terms with $s\ge1$ do
not contribute to the following calculation, and as in (3.7), we
\vs{-5pt}have
$$%
\begin{array}{ll}
&((d\tau_2+\tau_2
d)\g)_{\l_1,...,\l_q}(J_{E_{j_1,k_1}}^{n_1},...,J_{E_{j_q,k_q}}^{n_q})
\vs{2pt}\\
=\!\!\!\!&\dis\sum_{i=1}^q\g_{\l_1,...,\l_q}
(J_{E_{j_1,k_1}}^{n_1},...,J_{E_{j_{i-1},k_{i-1}}}^{n_{i-1}},
[J_{E_{j_i,k_i}}^{n_i}{\ssc\,}_{\l_i}J_h^0],
J_{E_{j_{i+1},k_{i+1}}}^{n_{i+1}},...,J_{E_{j_q,k_q}}^{n_q})
\vs{0pt}\\
=\!\!\!\!&\dis\sum_{i=1}^q(k_i-j_i)\g_{\l_1,...,\l_q}
(J_{E_{j_1,k_1}}^{n_1},...,J_{E_{j_q,k_q}}^{n_q}).
\end{array}
\vs{-9pt}\eqno(3.14)$$%
Thus as in Lemma 2.2, we obtain the following lemma.%
\par%
{\it Lemma 3.3:} {\it A $q$-cocycle in $\wt D^q(gc_N,\C)$ is
equivalent to a homogenous $q$-cocycle
$\g$ of degree zero \vs{-9pt}satisfying %
$$%
\g_{\l_1,...,\l_q}
(J_{E_{j_1,k_1}}^{n_1},...,J_{E_{j_q,k_q}}^{n_q})=0\mbox{ \ \ if \
\ }\sum_{i=1}^q(j_i-k_i)\ne0.
\vs{-15pt}\eqno(3.15)$$%
\ }\qed%
\def\eta{\Delta}
For a $q$-cochain $\g\in\wt C^q(gc_N,\C)$, we define a linear map
$\eta\g:gc_N^{\otimes q}\to\C[\l_1,...,\l_{q-1}]$ by
$$%
\eta\g(a_1\otimes\cdots\otimes
a_q)\!=\!\g_{\l_1,...,\l_q}(a_1,...,a_q)|_{\l_q=-\l_1-...-\l_{q-1}
}\!=\!\g_{\l_1,...,\l_{q-1},-\l_1-...-\l_{q-1}
}(a_1,...,a_q), \eqno(3.16)
$$%
for $a_1,...,a_q\in gc_N$ (we define $\eta\g=\g$ if $q=0$, and
define $\eta\g(a_1)=\g_{\l_1}(a_1)|_{\l_1=0}$ if $q=1$). Let
$C'^q(gc_N,\C)=\{\eta\g\,|\,\g\in\wt C^q(gc_N,\C)\}$. Then we
obtain a linear map $\eta:\wt C^q(gc_N,\C)\to C'^q(gc_N,\C)$. If
$\g\in\ptl\wt C^q(gc_N,\C)=(\sum_{i=1}^q\l_i)\wt C^q(gc_N,\C)$
(note that $\ptl_\C=0$, cf.~(2.15)), then clearly $\eta\g=0$. Thus
$\eta$ factors to a map $\eta: C^q(gc_N,\C)\to C'^q(gc_N,\C)$.
\par%
{\it Lemma 3.4: The map $\eta:C^q(gc_N,\C)\to C'^q(gc_N,\C)$ is an
isomorphism as spaces.\par%
Proof:} %
Suppose $\eta\g=0$ for a $q$-cochain $\g$. For $a_1,...,a_q\in
gc_N$, regarding $\g_{\l_1,...,\l_q}(a_1,...,a_q)$ as a polynomial
in $\l_q$, we see that it has a root $\l_q=-\sum_{i=1}^{q-1}\l_i$,
i.e., it is divided by $\sum_{i=1}^q\l_i$. \vs{-9pt}Thus
$$
\phi(a_1\otimes\cdots\otimes
a_q)=(\sum_{i=1}^q\l_i)^{-1}\g_{\l_1,...,\l_q}(a_1,...,a_q),
\vs{-9pt}\eqno(3.17)%
$$
defines a map $\phi:gc_N^{\otimes N}\to\C[\l_1,...,\l_q]$.
Obviously, $\phi$ is a $q$-cochain, and
$\g=(\sum_{i=1}^q\l_i)\phi\in\ptl\wt C^q(gc_N,\C)$.\hfill$\Box$.
\par%
Thus we can identify $C^q(gc_N,\C)$ with the space
$C'^q(gc_N,\C)$. We call an element in $C'^q(gc_N,\C)$ a {\it
reduced \mbox{$q$-cochain}}. We define the operator
$d:C'^q(gc_N,\C)\to C'^{q+1}(gc_N,\C)$ by $d\eta\g=\eta d\g$, and
then we have similar notions of {\it reduced $q$-cocycles, reduced
$q$-coboundaries}.%
\par%
{\it Lemma 3.5: Theorem 2.5(2) holds.\par%
Proof:} %
Clearly, by (2.14), $\wt D^0(gc_N,\C)=\wt
C^0(gc_N,\C)=\C$, and $\wt B^0(gc_N,\C)=0$. Thus $\wt
H^0(gc_N,\C)=\C$. Also by (2.15), $\ptl\wt C^0(gc_N,\C)=0$ and we
have $H^0(gc_N,\C)=\C$.
\par
Suppose $\g\in\wt C^1(gc_N,\C)$ such that $d\g\in\ptl\wt
C^2(gc_N,\C)$, i.e., there is $\phi\in\wt C^2(gc_N,\C)$ such
\vs{-4pt}that
$$%
\begin{array}{ll}
\g_{\l_1+\l_2}([u_{\l_1}v])\!\!\!\!&=
-(d\g)_{\l_1,\l_2}(u,v)=-(\ptl\phi)_{\l_1,\l_2}(u,v)\vs{4pt}\\ &%
=-
(\ptl_{\C}+\l_1+\l_2)\phi_{\l_1,\l_2}(u,v)=-(\l_1+\l_2)\phi_{\l_1,\l_2}(u,v),
\end{array}
\vs{-4pt}\eqno(3.18)$$%
(cf.~(2.13) and (2.15)) for $u,v\in S_N$. By (2.6), we
\vs{-9pt}have
$$%
[J_A^n{\ssc\,}_{\l_1} J^0]=\sum_{s=1}^n\biggl(\!\begin{array}{c}n\\
s\end{array}\!\biggr)(\l_1+\ptl)^sJ_A^{n-s},%
\vs{-9pt}\eqno(3.19)
$$%
for $A\in gl_N,\,n\in\Z_+$, where $J^0=J_I^0$.
\vs{-7pt}Thus by (2.11), (3.18) and (3.19), we have%
$$%
\sum_{s=1}^n\biggl(\!\begin{array}{c}n\\
s\end{array}\!\biggr)(-\l_2)^s\g_{\l_1+\l_2}(J_A^{n-s})=\g_{\l_1+\l_2}
[J_A^n{\ssc\,}_{\l_1}J^0]=-(\l_1+\l_2)\phi_{\l_1,\l_2}(u,v).
\vs{-7pt}\eqno(3.20)$$%
Let $\l_1=\l-\l_2$, then expressions in (3.20) are polynomials in
$\l,\l_2$ and the right-hand side is divided by $\l$, thus each
term in the left-hand side is divided by $\l$. Therefore we can
set %
$$\g'{\ssc\,}_\l(J_A^n)=\l^{-1}\g{\ssc\,}_\l(J_A^n)\mbox{ \  for \ }a\in gl_N,\,n\in\Z_+.\eqno(3.21)$$
Clearly, (3.21) defines a $1$-cochain $\g'\in\wt C^1(gc_N,\C)$,
and we have $\g=\ptl\g'\in\ptl\wt C^1(gc_N,\C)$. This proves that
$H^1(gc_N,\C)=0$.
\par%
Now suppose $\g\in\wt D^1(gc_N,\C)$ is a $1$-cocycle. This means
that $\phi=0$ in (3.18) and (3.20), and so, we obtain $\g=0$. Thus
$\wt H^1(gc_N,\C)=0$.
\par%
Next suppose $\psi\in\wt D^2(gc_N,\C)$ is a homogenous $2$-cocycle
of degree zero. We define a $1$-cochain $f$ which is uniquely
determined by
$$%
f_{\l_1}(J_A^n)=(n+1)^{-1}\ptla\psi_{\l_1,\l}(J_A^{n+1},J^0)|_{\l=0}.
\eqno(3.22)
$$%
Set $\g=\psi+df$, which is also a homogenous $2$-cocycle of degree
zero. Then
$$%
\ptla\g_{\l_1,\l}(J_A^n,J^0)|_{\l=0}=\ptla\psi_{\l_1,\l}(J_A^n,J^0)|_{\l=0}
-\ptla f_{\l_1+\l}([J_A^n{\ssc\,}_{\l_1}J^0])|_{\l=0}=0,
\eqno(3.23)$$%
where, the last equality follows from (3.19), (2.11) and (3.22) if
$n\ge1$, or from the fact that $\psi_{\l_1,\l}(J_A^0,J^0)$ is a
constant polynomial (cf.~(3.3)) if $n=0$. Thus we have%
\par\ni$%
\begin{array}{ll}%
0\!\!\!\!\!&=\!\ptla(d\g)_{\l_1,\l_2,\l}(J_A^m,J_B^n,J^0)|_{\l=0}\vs{4pt}\\
&=\!
\ptla(-\g_{\l_1+\l_2,\l}([J_A^m{\ssc\,}_{\l_1}J_B^n],J^0)\!+\!
\g_{\l_1+\l,\l_2}([J_A^m{\ssc\,}_{\l_1}J^0],J_B^n)
\!-\!\g_{\l_2+\l,\l_1}([J_B^n{\ssc\,}_{\l_2}J^0],J_A^m))|_{\l=0}\vs{4pt}\\
&=\!m\g_{\l_1,\l_2}(J_A^{m-1},J_B^n)+n\g_{\l_1,\l_2}(J_A^m,J_B^{n-1}),
\end{array}%
$\hfill(3.24)\par\ni%
for $A,B\in gl_N,\,m,n\in\Z_+$, where the second equality follows
from (2.13), the last equality follows from (3.23), (3.19) and
(2.11). Induction on $n\ge0$ in (3.24) proves
$\g_{\l_1,\l_2}(J_A^m,J_B^n)=0$. Thus $\g=0$ and so $\wt
H^2(gc_N,\C)=0$.
\par%
Finally, suppose $\psi'=\eta\psi\in C'^2(gc_N,\C)$ is a reduced
$2$-cochain. By (2.13)
and (3.16), 
$$%
(d\psi')_{\l_1,\l_2}(a_1,a_2,a_3)=
-\psi'_{\l_1+\l_2}([a_1{\ssc\,}_{\l_1}a_2],a_3)
+\psi'_{-\l_2}([a_1{\ssc\,}_{\l_1}a_3],a_2)
-\psi'_{-\l_1}([a_2{\ssc\,}_{\l_2}a_3],a_1),
\eqno(3.25)$$%
for $a_1,a_2,a_3\in gc_N$. We define a reduced $1$-cochain
$f'=\eta f\in C'^1(gc_N,\C)$ as follows (note from (3.16) that
$f'(a)=f_\l(a)|_{\l=0}=f_0(a)$ is simply a linear function
$f':gc_N\to\C$, and it is not necessary to write down  explicitly
its representative (basic) $1$-cochain $f$)
$$%
f'(J_A^m)=(m+1)^{-1}\mbox{$d\over d\l$}\psi'_\l(J_A^m,J)|_{\l=0},
\eqno(3.26)
$$%
(recall (3.6) that $J=J_I^1$) for $A\in gl_N,\,m\in\Z_+$ (note
from (2.11) that $f'(\ptl a)=f_0(\ptl
a)=0$). By (2.13) and (3.16), %
$$%
(df')_\l(a_1,a_2)=-f'([a_1{\ssc\,}_\l a_2]), \eqno(3.27)
$$%
for $a_1,a_2\in gc_N$.\par%
 Now suppose $\psi'$ is a reduced $2$-cocycle. Then
$\g'=\psi'+df'$ is a reduced $2$-cocycle equivalent to $\psi'$. By
(3.26), (3.27) and (3.6), $$%
\mbox{$d\over d\l$}\g'_\l(J_A^m,J)|_{\l=0}=0\mbox{ \ for \ }A\in
gl_N,\,m\in\Z_+.
\eqno(3.28)$$%
Thus by (3.25), $$%
\begin{array}{ll}%
0\!\!\!\!&=\ptla(d\g')_{\l_1,\l}(J_A^m,J_B^n,J)|_{\l=-\l_1}\vs{4pt}\\
&= \ptla(-\g'_{\l_1+\l}([J_A^m{\ssc\,}_{\l_1}J_B^n],J)
+\g'_{-\l}([J_A^m{\ssc\,}_{\l_1}J],J_B^n)
-\g'_{-\l_1}([J_B^n{\ssc\,}_\l
J],J_A^m))|_{\l=-\l_1}\vs{4pt}\\
&=\ptla((m(\l_1+\l)+\l_1)\g'_{-\l}(J_A^m,J_B^n)
-((n(\l+\l_1)+\l)\g'_{-\l_1}(J_B^n,J_A^m))|_{\l=-\l_1},
\end{array}%
\eqno(3.29)$$%
where the last equality follows from (3.28) and
(3.6) (similarly to the discussion after (3.7), $\l_1+\ptl$ and
$\l+\ptl$ can be replaced by $\l+\l_1$ and the terms with $s\ge2$
do not contribute to the calculation). Using (2.12) and (3.16),
the right-hand side
of (3.29) is equal to%
$$
(m+n+1)\g'_{\l_1}(J_A^m,J_B^n)-\l_1\mbox{$\ptl\over
\ptl\l_1$}\g'_{\l_1}(J_A^m,J_B^n)=0. \eqno(3.30)
$$
{}From (3.30), we obtain
$$%
\g'_\l(J_A^m,J_B^n)=c^{(m,n)}_{A,B}\l^{m+n+1} \mbox{ \ for some \
}c^{(m,n)}_{A,B}\in\C.%
\eqno(3.31)
$$%
In particular,
$$%
\mbox{$d\over d\l$}\g'_\l(J_A^m,J^0)|_{\l=0}=\d_{m,0}c_A,
\eqno(3.32)$$%
where $c_A=c_{A,I}^{(0,0)}$. Similarly to (3.29) (also cf.~(3.24)),
$$%
\begin{array}{ll}%
0\!\!\!\!&=\ptla(d\g')_{\l_1,\l}(J_A^m,J_B^n,J^0)|_{\l=-\l_1}\vs{4pt}\\
&=
-\ptla\g'_{\l_1+\l}([J_A^m{\ssc\,}_{\l_1}J_B^n],J^0)|_{\l=-\l_1}
+m\g'_{\l_1}(J_A^{m-1},J_B^n)+n\g'_{\l_1}(J_A^m,J_B^{n-1})\vs{4pt}\\
&=-(^{\ \,m}_{m+n})\l_1^{m+n}c_{AB}+ (^{\ \
n}_{m+n})(-\l_1)^{m+n}c_{BA}
+(mc_{A,B}^{(m-1,n)}+nc_{A,B}^{(m,n-1)})\l_1^{m+n},
\end{array}%
\eqno(3.33)$$%
where the last equality follows from (2.6), (2.11), (3.16), (3.31)
and (3.32). Taking $m=n=0$, we obtain $c_{AB}=c_{BA}$. Thus
$$%
mc_{A,B}^{(m-1,n)}+nc_{A,B}^{(m,n-1)}=((^{\
\,m}_{m+n})-(-1)^{m+n}(^{\ \ n}_{m+n}))c_{AB}. \eqno(3.34)
\vs{-4pt}$$%
Thus we \vs{-4pt}solve
$$
c_{A,B}^{(m,n)}=(-1)^n{m!n!\over(m+n+1)!}c_{AB}\mbox{ \ for \
}A,B\in gl_N,\,m,n\in\Z_+. %
\vs{-4pt}\eqno(3.35)$$%
{}From (3.31) and the fact that $c_A=c_{A,I}^{(0,0)}$ and that
$c_{AB}=c_{BA}$, we see that the map $A\mapsto c_A$ is a trace of
$gl_N$, i.e., $c_A$ is a scalar multiple of ${\rm tr}(A)$ for $\in
gl_N$. Thus (3.31) and (3.35) show that $\g'$ is a multiple of
$\psi'$ which is defined \vs{-4pt}by
$$
\psi'_\l(J_A^m,J_B^n)=(-1)^n{m!n!\over(m+n+1)!}{\rm
tr}(AB)\l^{m+n+1}. \vs{-4pt}\eqno(3.36)
$$
To see that $\psi'$ is a nontrivial reduced $2$-cocycle, first
\vs{-4pt}define
$$
\psi_{\l_1,\l_2}(J_A^m,J_B^n)=(-1)^n{m!n!\over(m+n+1)!}
((-1)^m\l_1^{m+n+1}-(-1)^n\l_2^{m+n+1}){\rm tr}(AB)\l^{m+n+1}.
\vs{-4pt}\eqno(3.37)$$%
Clearly, $\psi$ is a $2$-cochain (recall the second sentence in
the paragraph before (3.2)), and $\psi'=\eta\psi$ is a reduced
$2$-cochain. One can easily check that $d\psi'=0$ and that
$\psi'\ne df'$ for any reduced $1$-cochain $f'$. This proves that
$H^2(gc_N,\C)=\C\psi'$.\qed%
{\it Lemma 3.6: Theorem 2.5(3) holds.\par%
Proof: }%
We define an operator $\tau:\wt C^q(gc_N,\C_a)\to\wt
C^{q-1}(gc_N,\C_a)$ by
$$%
(\tau\g)_{\l_1,...,\l_{q-1}}(a_1,...,a_{q-1})=(-1)^{q-1}\g_{\l_1,...,\l_{q-1},\l}
(a_1,...,a_{q-1},J)|_{\l=0}, \eqno(3.38)
$$%
for $a_1,...,a_{q-1}\in gc_N$. Similarly to the discussions in (3.7)
and (3.8), we \vs{-7pt}have
$$%
\begin{array}{ll}%
 ((d\tau+\tau
d)\g)_{\l_1,...,\l_q}(J_{A_1}^{n_1},...,J_{A_q}^{n_q})\!\!\!\!&\dis=(\sum_{i=1}^q\l_q)
\g_{\l_1,...,\l_q}(J_{A_1}^{n_1},...,J_{A_q}^{n_q})\\
&\equiv -a\g_{\l_1,...,\l_q}(J_{A_1}^{n_1},...,J_{A_q}^{n_q}) \
({\rm mod\,}\ptl\wt C^q(gc_N,\C_a)\,),
\end{array}%
 \eqno(3.39)$$%
(note that $\ptl\wt C^q(gc_N,\C_a)=(a+\sum_{i=1}^q\l_i)\wt
C^q(gc_N,\C_a)$ by (2.15), since $\ptl_{\C_a}=a\ssc\,$). Now
suppose $\g\in\wt C^q(gc_N,\C_a)$ such that $d\g\in\ptl\wt
C^{q+1}(gc_N,\C_a)$, i.e., there exists a $(q+1)$-cochain $\phi$
such that $d\g=(a+\sum_{i=1}^{q+1}\l_i)\phi$. Clearly, by (3.38)
$\tau d\g=(a+\sum_{i=1}^q\l_i)\tau\phi\in\ptl\wt C^q(gc_N,\C_a)$.
Thus (3.39) shows that $\g\equiv -d(a^{-1}\tau\g)\ ({\rm
mod\,}\ptl\wt C^q(gc_N,\C_a)\,)$ is a reduced coboundary (note
that we assume $a\ne0\ssc\,$),
i.e., $H^q(gc_N,\C_a)=0$.\qed%
{\it Lemma 3.7: Theorem 2.5(4) holds.%
\par Proof: }
Note that as spaces, we have
$\C_a^N[\ptl][\l_1,...,\l_q]=\C^N[\l_1,...,\l_q,\ptl]$, and  a
$q$-cochain $\wt\g\in\wt C^q(gc_N,\C_\a^N[\ptl])$ can be regarded
as a map $\wt\g:gc_N^{\otimes q}\to \C^N[\l_1,...,\l_q,\ptl]$,
$$
\wt\g(a_1\otimes\cdots\otimes
a_q)=\wt\g_{\l_1,...,\l_q,\ptl}(a_1,...,a_q), \eqno(3.40)
$$%
for $a_1,...,a_q\in gc_N$. Regarding (3.40) as a polynomial in
$\l_1,...,\l_q,\ptl$ with coefficients in $\C^N$, then similarly to
Lemma 3.4, a reduced\vs{-7pt} $q$-cochain%
 $$\g\in
C^q(gc_N,\C_\a^N[\ptl])=\wt
C^q(gc_N,\C_\a^N[\ptl])/(\ptl+\sum_{i=1}^q\l_i)\wt
C^q(gc_N,\C_\a^N[\ptl]),%
\vs{-7pt}\eqno(3.41)$$%
is uniquely determined by the coefficient of $\ptl^0$ in (3.40).
Thus a reduced $q$-cochain $\g$ can be regarded as a map
 $\g:gc_N^{\otimes q}\to \C^N[\l_1,...,\l_q]$,
$$
\g(a_1\otimes\cdots\otimes a_q)=\g_{\l_1,...,\l_q}(a_1,...,a_q).
\eqno(3.42)
$$
Define an operator $\tau_0:C^q(gc_N,\C_\a^N[\ptl])\to
C^{q-1}(gc_N,\C_\a^N[\ptl])$ by (cf.~(3.38))
$$%
(\tau_0\g)_{\l_1,...,\l_{q-1}}(a_1,...,a_{q-1})=(-1)^{q-1}
\g_{\l_1,...,\l_{q-1},\l}(a_1,...,
a_{q-1},J^0)|_{\l=0}.\eqno(3.43)
$$%
Similarly to the discussions in (3.7) and (3.8), using (3.19), we
have (comparing with (3.7), all terms corresponding to the
right-hand side of (3.7) are now zero because $J$ has been
replaced by $J^0$ and we do not take partial derivative $\ptla$;
but note that since the first sum in (2.13) is not zero in this
case, we have one more
term here)%
$$%
((d\tau_0+\tau_0d)\g)_{\l_1,...,\l_q}(J_{A_1}^{n_1},...,J_{A_q}^{n_q})=
J^0{\ssc\,}_\l\g_{\l_1,...,\l_q}(J_{A_1}^{n_1},...,J_{A_q}^{n_q})|_{\l=0}.
\eqno(3.44)$$%
Now by (2.9), the $\l$-action of $gc_N$ on its module
$\C_\a^N[\ptl]$ in particular satisfies $J^0{\ssc\,}_\l v=v$ for
$v\in\C^N_\a[\ptl]$. Thus the right-hand side of (3.44) is simply
$-\g_{\l_1,...,\l_q}(J_{A_1}^{n_1},...,J_{A_q}^{n_q})$, i.e., we
\vs{-5pt}obtain
$$
\g=(d\tau_0+\tau_0d)\g.\vs{-5pt}\eqno(3.45)
$$
In particular, if $\g$ is a reduced cocycle, (3.45) gives that
$\g=d(\tau_0\g)$ is a coboundary, i.e.,
$H^q(gc_N,\C_\a^N[\ptl])=0$.\par%
Clearly, the above proof works for any $gc_N{\ssc\,}$-module $M$
satisfying the condition stated
in Theorem 2.5(4).\qed%
Thus Theorem 2.5(2)-(4) is proved. \par%
\vskip-5pt\ \par\ni{\bf IV. \ PROOF OF THEOREM 2.5(1)}\par%
This section is devoted to the proof of Theorem 2.5(1). By Lemma
3.5, it remains to consider the case $q=3$. Since some of the
following arguments also work for general \mbox{$q$-cocycles}, we
shall first consider $q$-cocycles with $q\ge3$ so that it may be
possible to use
these arguments to determine higher dimensional cohomologies in the future. %
\par%
Let $gc=gc_1$. It has a free generating set
$S=\{J^n\,|\,n\in\Z_+\},$ such \vs{-7pt}that
$$%
[J^m{\ssc\,}_\l
J^n]=\sum_{s=1}^m\biggl(\!\begin{array}{c}m\\
s\end{array}\!\biggr)(\l+\ptl)^sJ^{m+n-s}-\sum_{s=1}^n\biggl(\!\begin{array}{c}n\\
s\end{array}\!\biggr)(-\l)^s
J^{m+n-s},%
\vs{-7pt}\eqno(4.1)$$%
for $m,n\in\Z_+$. We shall give some more notations. An element in
$\Z_+^q$ is denoted by
$$%
\nn=\nn[q]=(n_1,...,n_q),\;\;\;\;n_1,...,n_q\in\Z_+, \eqno(4.2)
$$%
(when there is no confusion we denote it by $\nn$, otherwise we
denote it by $\nn[q]$). Denote $J^\nn=J^{n_1}\otimes\cdots\otimes
J^{n_q}=(J^{n_1},...,J^{n_q})\in gc^{\otimes q}$. Denote
$\ll=\ll[q]=(\l_1,...,\l_q)$. For $\nn\in\Z_+^q$, let $|\ul
n|=\sum_{i=1}^qn_i$, called the {\it level of $\,\ul n$}.
We define a total ordering on $\Z_+^q$ by the {\it level-lexicographical order}, i.e.,%
$$%
\mm\!<\!\nn\ \Lra\ |\mm|\!<\!|\nn|,\mbox{ or
}|\mm|\!=\!|\nn|\mbox{ and }\exists\,p\mbox{ such
that }m_i\!=\!n_i\mbox{ for }i\!<\!p\mbox{ and }m_p\!<\!n_p,%
\eqno(4.3)
$$%
for $\ul m,\ul n\in\Z_+^q$. Set
$$\NN=\{\ul n\in\Z_+^q\,|\,n_1\le n_2\le...\le
n_q\}.\eqno(4.4)$$%
For $m,n\in\Z$, we denote $[m,n]=\{m,m+1,...,n\}$. Let $\SS$ be
the permutation group on the index set $[1,q]$, which acts on
$\C^q$ by $\si(v)=(v_{\si(1)},...,v_{\si(q)})$ for
$v=(v_1,...,v_q)\in\C^q$. Then for any $\nn\in\Z_+^q$, there is a
unique $\nn^*\in\NN$ and some $\si\in\SS$ such that
$\nn^*=\si(\nn)\in\NN$ and $\nn^*\le\si(\nn)$. In fact
$$%
\nn^*={\rm min}\{\si(\nn)\,|\,\si\in\SS\},
\eqno(4.5)$$%
is the minimal element in $\SS(\nn)=\{\si(\nn)\,|\,\si\in\SS\}$.
\par%
A $q$-cochain $\g$ is uniquely determined by $\g_\ll(J^{\ul n})$
for $\ul n\in\NN$ \vs{-4pt}and
$$%
\g_\ll(J^{\ul n})={\rm sgn}(\si)\g_{\si(\ll)}(J^{\si(\ul n)}),
\vs{-4pt}\eqno(4.6)
$$%
for $\ul n\in\Z_+^q,\,\si\in\SS$, where ${\rm sgn}(\si)$ is the
signature of the permutation $\si$. In fact, $\g_\ll(J^\nn)$ can
be arbitrary polynomial in $\ll$ satisfying (4.6) for all $\si$
such that $\si(\nn)=\nn$.
\par%
First we construct a $3$-cochain $\bar\g$ as \vs{-7pt}follows:
$$%
\bar\g_\ll(J^{\ul n})=
\left\{\begin{array}{ll}%
\l_2^{n_3}-\l_1^{n_3}&\mbox{if \ }n_1=n_2=0,n_3\ne0,\\
0&\mbox{otherwise},
\end{array}\right.%
\vs{-7pt}\eqno(4.7)$$%
for $\ul n\in{\cal N}_3$ (note that in the first case, we let
$n_3\ne0$ in order to avoid the problem on how to deal with $0^0$
when we set $\l_1=0$).
\par%
{\it Lemma 4.1: $\bar\g$ is a nontrivial $3$-cocycle.\par%
Proof: }%
One can define a {\it Leibniz $q$-cochain} by removing the
skew-symmetric condition (2.12), and define the {\it Leibniz
differential operator $d_L$} by changing (2.13) into
\par\ni$%
\begin{array}{ll}%
&\dis(d_L\g)_{\l_1,...,\l_{q+1}}(a_1,...,a_{q+1})\vs{2pt}\\
=\!\!\!\!&\dis
 \sum_{i=1}^{q+1}(-1)^{i+1}a_i{\ssc\,}_{\l_i}
\g_{\l_1,...,\hat\l_i,...,\l_{q+1}}(a_1,...,\hat
a_i,...,a_{q+1})\vs{0pt}\\ %
&+\dis\sum_{1\le i<j\le
q+1}(-1)^i\g_{\l_1,...,\hat\l_i,...,\l_{j-1},\l_i+\l_j,\l_{j+1},...,\l_{q+1}}
(a_1,...,\hat
a_i,...,a_{j-1},[a_i{\ssc\,}_{\l_i}a_j],a_{j+1},...,a_{q+1}),
\end{array}%
\!\!\!\!\!\!$\hfill(4.8)\par\ni%
(note that if $\g$ is a (regular) $q$-cochain, then (4.8)
coincides with (2.13), i.e., $d=d_L$ in this case). Then we obtain
{\it Leibniz cohomology}  (cf.~Ref.~2). We shall not discuss
Leibniz cohomology here, but we define a Leibniz $2$-cochain $f$
by
$$
f_{\l_1,\l_2}(J^0,J^0)=1\mbox{ \ and \
}f_{\l_1,\l_2}(J^m,J^n)=0\mbox{ \ if \ }(m,n)\ne(0,0).
\eqno(4.9)$$%
One can immediately check that $\bar\g=d_Lf$ (thus $\bar\g$ is a
Leibniz $3$-coboundary). Therefore $d\bar\g=dd_Lf=d_L^2f=0$, i.e.,
$\bar\g$ is a (regular) $3$-cocycle. However there is no
$2$-cochain $\phi$ such that $d\phi=\bar\g$ because if
$d\phi=\bar\g$ then we also have $\phi_{\l_1,\l_2}(J^0,J^0)=1$ and
so $\phi$ is not a (regular) $2$-cochain ((2.12) is not
satisfied). Thus $\bar\g$ is a nontrivial
$3$-cocycle.\qed%
Now let $\g$ be a $q$-cocycle with $q\ge3$. By Lemma 3.2, we can
suppose $\g$ is homogenous with degree zero. First
we have the following lemma.\par%
{\it Lemma 4.2: If $q=3$, by replacing $\g$ by $\g-c\bar\g$ for
some $c\in\C$,
we can suppose $\g_\ll(J^0,J^0,J)=0.$\par%
Proof: }%
Note that $\g_\ll(J^0,J^0,J)$ is a linear polynomial in $\ll$
(cf.~(3.3)) which is skew-symmetric with respect to $\l_1,\l_2$ by
(2.12). Thus $\g_\ll(J^0,J^0,J)=c(\l_2-\l_1)$ for some $c\in\C$.
Replacing $\g$ by $\g-c\bar\g$, we have the lemma.\qed%
To prove (2.16), our strategy is the following: We want to prove
by induction on $\ul n\in\NN$ (with respect to the order (4.3))
that after a number of steps in each of which $\g$ is replaced by
$\g-\g'$ for some $q$-coboundaries $\g'$ we obtain that
$\g_\ll(J^{\ul m})=0$ for all $\ul m\in\NN,\,\ul m\le\ul n$ (thus
we obtain that $\g_\ll(J^{\ul n})=0$ for all $\ul n\in\NN$, i.e.,
$\g=0$, after a countably infinite number of steps; this amounts
to saying that $\g$ is subtracted by an infinite sum of
\mbox{$q$-coboundaries}, but from the following proof we see that
this infinite sum is summable, cf.~the statement after (3.4)). For
the case $q=3$, this will be done by a number of lemmas
(unfortunately, not all arguments work for $q\ge4$, cf.~the proof
of Lemma 4.6).
\par%
{\it Lemma 4.3: $\g_\ll(J^\nn)=0$ if $|\nn|\le1$. \par%
Proof: }%
Note that $\g_\ll(J^\nn)$ is a polynomial in $\ll$ on degree
$|\nn|$. If $|\nn|=0$, we have $\g_\ll(J^\nn)=0$ by (4.6). If
$|\nn|=1$, then $\nn=(0,...,0,1)$ and $\g_\ll(J^\nn)$ is
skew-symmetric with respect to $\l_1,...,\l_{q-1}$, thus divided
by $\prod_{1\le i<j\le q-1}(\l_i\!-\!\l_j)$, which has degree
$(q\!-\!1)(q\!-\!2)/2>1$ if $q\!>\!3$. Thus $\g_\ll(J^\nn)=0$ if
$q>3$. If $q=3$, then $\g_\ll(J^\nn)=\g_\ll(J^0,J^0,J)=0$ by Lemma
4.2.\qed%
Now suppose $|\nn|\ge2$. We set
$i_0=\#\{i\in[1,q]\,|\,n_i=0\}\ge0$ (where $\#X$ stands for the
size of the finite set $X$), $i_2=q-\#\{i\in[1,q]\,|\,n_i=n_q\}\le
q-1.$ If $i_0\ne i_2$, we set $i_1$ to satisfy
$$%
0=n_1=...=n_{i_0}<n_{i_0+1}\le...\le
n_{i_1}<n_{i_1+1}=...=n_{i_2}<n_{i_2+1}=...=n_q;%
\eqno(4.10)$$%
if $i_0=i_2$, we set $i_1=0$.
\par%
Let $\mm\in{\cal N}_{q+1}$ be such that $|\mm|=|\nn|+1$. Consider
$(d\g)_{\ll[q+1]}(J^\mm)$ (cf.~notation (4.2)). Note that when we
substitute (4.1) into (2.13), using (2.11) and (2.12), we obtain
that $(d\g)_{\ll[q+1]}(J^\mm)$ is a combination of
$\g_{\ll'}(J^\kk)$ with coefficients being polynomials in
$\ll[q\!+\!1]$, where $\kk\!\in\!\NN,\,|\kk|\!\le\!|\nn|$, and
$\ll'\!=\!(\l'_1,...,\l'_q)$ such that each $\l_i'$ is a linear
polynomial in $\ll[q\!+\!1]$. Using the inductive assumption,
$\g_{\ll'}(J^\kk)=0$ if $|\kk|<|\nn|$. Thus the terms with $s\ge2$
in (4.1) do not contribute to (2.13) (cf.~the discussion after
(3.7)), and so we have (here we use (4.8) instead of
\vs{-3pt}(2.13))
$$%
\begin{array}{rl}%
0=\!\!\!\!&\dis(d\g)_{\ll[q+1]}(J^\mm)\vs{4pt}\\
=\!\!\!\!&\dis\sum_{1\le i<j\le q+1}(-1)^i(m_j\l_i-m_i\l_j)\
\g_{\l_1,...,\hat\l_i,...,\l_{j-1},\l_i+\l_j,\l_{j+1},...,\l_{q+1}}
(J^{\mm(i,j)}),
\vs{4pt}\\ %
\end{array}%
\vs{-8pt}\eqno(4.11)
$$%
where $$\mm(i,j)=(m_1,...,\hat
m_i,...,m_{j-1},m_i+m_j-1,m_{j+1},...,m_{q+1}),%
\eqno(4.12)$$%
and the right-hand side of (4.11) is a combination of
$\g_{\ll'}(J^{\mm(i,j)^*})$ (cf.~(4.5) and (4.6)).
\par%
{\it Lemma 4.4: $\g_\ll(J^\nn)=0$ if $n_1\ge1$ (i.e., $i_0=0$).
\par%
Proof: }%
In (4.11), take $\mm=(0,n_1,...,n_{q-1},n_q+1)\in{\cal N}_{q+1}$.
In (4.12), if $i\ne1$, then $m_1=0<n_1$ and so
$\mm(i,j)^*\le\mm(i,j)<\nn$; by induction,
$\g_{\ll'}(J^{\mm(i,j)^*})=0$. Similarly,
$\g_{\ll'}(J^{\mm(i,j)^*})=0$ if $j\ne q+1$. Thus the only
possible nonzero term in (4.11) is the one with $(i,j)=(1,q+1)$.
Since $\mm(1,q+1)=\nn$, (4.11) \vs{-7pt}gives
$$-(n_q+1)\l_1\g_{\l_2,...,\l_q,\l_1+\l_{q+1}}(J^\nn)=0.\vs{-7pt}\eqno(4.13)$$
This gives
the lemma.\qed%
From now on, we assume that $n_1=0$.
\par%
{\it Lemma 4.5: $\g_\ll(J^\nn)=0$ if $\,i_2=q-1$ and $n_q\ge n_{q-1}+2$ (cf.~(4.10)).\par%
Proof: }%
As above, now (4.11) \vs{-9pt}gives (cf.~(4.13))
$$%
\sum_{i=1}^{i_0+1}(-1)^i(n_q+1)\l_i
\g_{\l_1,...,\hat\l_i,...,\l_q,\l_i+\l_{q+1}}(J^\nn)=0.
\vs{-8pt}\eqno(4.14)$$%
Replacing $(\l_1,...,\l_{q+1})$ by $(\l,\l_1,...,\l_q)$ and
applying the operator $\ptla|_{\l=0}$ to (4.14),
we \vs{-8pt}obtain%
$$%
\g_\ll(J^\nn)=-\sum_{i=1}^{i_0}(-1)^i\l_i\ptla\g_{\l,\l_1,...,\hat\l_i,...,\l_{q-1},
\l_i+\l_q}(J^\nn)|_{\l=0}. \vs{-8pt}\eqno(4.15)
$$%
We define a $(q-1)$-cochain $f$ as \vs{-7pt}follows%
$$%
f_{\ll[q-1]}(J^\kk)=\left\{\begin{array}{ll}%
-n_q^{-1}\ptla\g_{\l,\l_1,...,\l_{q-1}}(J^\nn)|_{\l=0}&\mbox{
if \ \ }\kk=\nn^-,\vs{4pt}\\
0&\mbox{ otherwise,}
\end{array}\right.%
\vs{-7pt}\eqno(4.16)$$%
for $\kk\!\in\!{\cal N}_{q-1}$, where
$\nn^-\!=\!(n_2,...,n_{q-1},n_q-1)\!\in\!{\cal N}_{q-1}$
(cf.~(4.10)). Indeed, $f$ is a \mbox{$(q\!-\!1)$-cochain} (cf.~the
statement after (4.6)): Write $\nn^-$ as $\nn^-=(n_1^-,...,n_{q-1}^-)$, then
$n_i^-=n_j^-\ \Lra\ n_{i+1}=n_{i+1}$. Thus the skew-symmetric
condition (4.6) for $f$ follows from the skew-symmetric condition
for $\g$. We claim that
$$%
\g_\ll(J^\kk)=(df)_\ll(J^\kk),\eqno(4.17)
$$%
for all $\kk\in\NN$ with $\kk\le\nn$. If $\kk=\nn$, similarly to
(4.14), we \vs{-8pt}have%
$$%
\begin{array}{ll}
(df)_\ll(J^\nn)\!\!\!\!&=
\dis\sum_{i=1}^{i_0}(-1)^in_q\l_if_{\l_1,...,\hat\l_i,...,
\l_{q-1},\l_i+\l_q}(J^{\nn^-})=\g_\ll(J^\nn),%
\end{array} \vs{-8pt}\eqno(4.18)$$%
where the last equality follows from (4.15) and (4.16). If
$\kk<\nn$, when we substitute (4.1) into (2.13) for
$(df)_\ll(J^\kk)$, as in (4.11), $(df)_\ll(J^\kk)$ is a
combination of the form $f_{\ll'}(J^{\kk(i,j)})$, and we see that
$\kk(i,j)<\nn^-$ (cf.~(4.12)), i.e., the term
$f_{\ll'}(J^{\nn^-})$ does not appear in $(df)_\ll(J^\kk)$, thus
$(df)_\ll(J^\kk)=0$, which is the same as $\g_\ll(J^\kk)$ by
inductive assumption. This proves (4.17). Thus by replacing $\g$
by $\g-df$, we have the lemma.\qed%
{\it Lemma 4.6: $\g_\ll(J^\nn)=0$ if $\,q=3$.\par%
Proof: }%
When $q=3$, by Lemmas 4.3-5, we are left to consider the cases
$\nn=(0,n_2,n_2)$ and $\nn=(0,n_2,n_2+1)$ for $n_2\ge1.$ First
suppose $\nn=(0,n_2,n_2)$. As in (4.14), we have
$$%
0=(d\g)_{\ll[4]}(J^0,J^0,J^{n_2},J^{n_2+1})
=(n_2+1)(-\l_1\g_{\l_2,\l_3,\l_1+\l_4}(J^\nn)+\l_2\g_{\l_1,\l_3,\l_2+\l_4}(J^\nn)).
\eqno(4.19)$$%
Setting $\l_4=0$, it gives that $\g_\ll(J^\nn)$ can be divided by
$\l_1$. So we can write $\g_\ll(J^\nn)=\l_1\g'_{\ll}$ for some
polynomial $\g'_\ll$, and (4.19) shows that
$\g'_{\l_2,\l_3,\l_1+\l_4}=\g'_{\l_1,\l_3,\l_2+\l_4}$. Setting
$\l_1=0$, this gives that
$\g'_{\l_2,\l_3,\l_4}=\g'_{0,\l_3,\l_2+\l_4}$. \vs{-7pt}Thus
$$\g_\ll(J^\nn)=\l_1\g'_{0,\l_2,\l_1+\l_3}.\vs{-7pt}\eqno(4.20)$$
But $\g_\ll(J^\nn)$ is skew-symmetric with respect to $\l_2,\l_3$,
we obtain $\g'_{0,\l_2,\l_1+\l_3}=-\g'_{0,\l_3,\l_1+\l_2}$.
Setting $\l_1=0$ and $\l_3=0$ respectively, we obtain that
$\g'_{0,\l_2,\l_3}=-\g'_{0,\l_3,\l_2}$ and
$\g'_{0,\l_2,\l_1}=-\g'_{0,0,\l_1+\l_2}$, which gives that
$\g'_\ll=0$. Thus $\g_\ll(J^\nn)=0$.
\par%
Next suppose $\nn=(0,n_2,n_2+1)$. We still have (4.20) for some
polynomial $\g'_\ll$. We assume that $n_2\ge2$ (the proof for the
case $n_2=1$ is similar and we leave it to the reader). For $1\le
i<n_2$, by (2.13) and the inductive assumption, we have
\par\ni\hs{-0.5ex}$
\begin{array}{ll}
0\!\!\!\!\!&=\!(d\g)_{\ll[4]}(J^0,J,J^{n_2-i},J^{n_2+i+1}) \vs{4pt}\\
&=\!(n_2\!-\!i)\l_1\g_{\l_2,\l_1+\l_3,\l_4}(J,J^{n_2-i-1},J^{n_2+i+1})
\!+\!(n_2\!+\!i\!+\!1)\l_1\g_{\l_2,\l_3,\l_1+\l_4}(J,J^{n_2-i},J^{n_2+i}).\end{array}
\!\!\!$\hfill(4.21) \par \ni%
Note that when $i=n_2-1$, the first term of the right-hand side is
zero since
$\g_\ll(J,J^0,J^{2n_2})=-\g_{\l_2,\l_1,\l_3}(J^0,J,J^{2n_2})$ and
$(0,1,2n_2)<\nn$. Thus induction on $i$ gives that
$\g_\ll(J,J^{n_2-i},J^{n_2+i})=0$. Then by (2.13) and the
inductive \vs{-5pt}assumption,
$$%
\begin{array}{rl}
0=\!\!\!\!&(d\g)_{\ll[4]}(J^0,J,J^{n_2},J^{n_2+1})\vs{4pt}\\
=\!\!\!\!&
-\l_1\g_{\l_1+\l_2,\l_3,\l_4}(J^\nn)-(n_2+1)\l_1\g_{\l_2,\l_3,\l_1+\l_4}
(J,J^{n_2},J^{n_2})\vs{4pt}\\ &
+(n_2\l_2-\l_3)\g_{\l_1,\l_2+\l_3,\l_4}(J^\nn)+((n_2+1)\l_2-\l_4)
\g_{\l_1,\l_3,\l_2+\l_4}(J^\nn).
\end{array}\vs{-5pt}
\eqno(4.22)$$%
Substituting (4.20) into (4.22), cancelling the common factor
$\l_1$, then setting $\l_1=\l_2=0$, we obtain that
$0=-(n_2+1)\g_{0,\l_3,\l_4}(J,J^{n_2},J^{n_2})-(\l_3+\l_4)\g'_{0,\l_3,\l_4}$,
which shows that $\g'_{0,\l_3,\l_4}$ is skew-symmetric with
respect to $\l_3,\l_4$. \vs{-5pt}Thus
$$%
f_{\l_1,\l_2}(J^{m_1},J^{m_2})=\left\{\begin{array}{ll}
-n_2^{-1}\g'_{0,\l_1,\l_2}&\mbox{if \ }(m_1,m_2)=(n_2,n_2),\vs{4pt}\\
0&\mbox{otherwise},
\end{array}\right.
\vs{-5pt}\eqno(4.23)$$%
defines a 2-cochain $f$ (cf.~(4.16)). Now as in the proof of Lemma
4.5, by replacing $\g$ by $\g-df$, we have the lemma. This also proves (2.16).\qed%
{\it Lemma 4.7: (2.17) holds.\par%
Proof: }%
By Lemma 3.5, it remains to consider the case $q=3$. Let $\bar\g$
be the 3-cocycle defined in (4.7). Let $\bar\g'=\eta\bar\g$ be the
corresponding reduced 3-cocycle (cf.~(3.16)). Clearly $\bar\g'$ is
nontrivial. Now suppose $\g'$ is arbitrary reduced 3-cocycle. As
in the paragraph before Lemma 4.3, we shall prove by induction on
$\nn\in{\cal N}_3$ that by replacing $\g'$ by
$\g'-c{\ssc\,}\bar\g'-df'$ for some $c\in\C$ and some reduced
$2$-cochain $f'$ we have $\g'_{\l_1,\l_2}(J^\mm)=0$ for
$\mm\le\nn$. Assume that we have proved
$\g'_{\l_1,\l_2}(J^\mm)=0$ for $\mm<\nn$.\par%
First suppose $\nn=(0,0,n_3)$. By (2.13), (3.16) and the inductive
assumption, we have
\par\ni\hs{1ex}%
$0\!=\!(d\g')_{\ll[3]}(J^0,J^0,J^0,J^{n_3+1})\!=\!
(n_3\!+\!1)(-\l_1\g'_{\l_2,\l_3}(J^\nn)
\!+\!\l_2\g'_{\l_1,\l_3}(J^\nn)\!-\!\l_3\g'_{\l_1,\l_2}(J^\nn)).
$\hfill(4.24)\par\ni%
\vs{-5pt}Thus
$$\g'_{\l_1,\l_2}(J^\nn)=\l_1\g'_{1,\l_2}(J^\nn)-\l_2\g'_{1,\l_1}(J^\nn).%
\vs{-5pt}\eqno(4.25)$$%
If $n_3=0$, then by (2.12) and (3.16),
$\g'_{\l_1,\l_2}(J^\nn)=-\g'_{\l_1,-\l_1-\l_2}(J^\nn)$, this
together with (4.25) gives that $\g'_{\l_1,\l_2}(J^\nn)=0$. If
$n_3=1$, then by (2.13), (3.16) and the inductive
\vs{-5pt}assumption,
$$%
\begin{array}{ll}0\!\!\!\!&=(d\g')_{\ll[3]}(J^0,J^0,J,J)\vs{4pt}\\
& =-\l_1(\g'_{\l_2,\l_1+\l_3}(J^\nn)
-\g'_{\l_2,-\l_2-\l_3}(J^\nn))+\l_2(\g'_{\l_1,\l_2+\l_3}(J^\nn)
-\g'_{\l_1,-\l_1-\l_3}(J^\nn)).
\end{array}\vs{-5pt}\eqno(4.26)$$%
Using (4.25) in (4.26), we see that $\g'_{1,\l_1}$ can be divided by
$\l_1$. Writing $\g'_{1,\l_1}=\l_1p(\l_1)$ for some polynomial
$p(\l_1)$ and using this in (4.26), cancelling the common factor
$\l_1\l_2$, and setting $\l_2=\l_3=0$, we see that
$p(\l_1)=c\in\C$ is a constant. Thus
$\g'_{\l_1,\l_2}(J^\nn)=c(\l_1-\l_2)$. Replacing $\g'$ by
$\g'+c\bar\g'$, we obtain that $\g'_{\l_1,\l_2}(J^\mm)=0$ for
$\mm\le\nn$.\par%
If $n_3\ge2$, we defines a reduced $2$-cochain $f'$ as
\vs{-5pt}follows:
$$%
f'_{\l_1}(J^{m_1},J^{m_2})=\left\{\begin{array}{ll} %
-\g'_{1,\l_1}(J^\nn)&\mbox{if \ }(m_1,m_2)=(0,n_3-1), \vs{4pt}\\
0&\mbox{otherwise},
\end{array}\right.
\vs{-5pt}\eqno(4.27)$$%
for $(m_1,m_2)\in{\cal N}_2$. Clearly, this indeed defines a
reduced $2$-cochain $f'$. Using (4.25), by replacing $\g'$ by
$\g'-df'$ as in the proof of Lemma 4.5, we have
$\g'_{\l_1,\l_2}(J^\mm)=0$ for $\mm\le\nn$.\par%
Next suppose $\nn=(0,n_2,n_2)$ for $n_2\ge1$. As in (4.25), from
$(d\g')_{\ll[3]}(J^0,J^0,J^{n_2},J^{n_2+1})=0$ we obtain that
$\g'_{\l_1,\l_2}(J^\nn)=\l_1\g'_{1,\l_2}(J^\nn)$. But
$\g'_{\l_1,\l_2}(J^\nn)=-\g'_{\l_1,-\l_1-\l_2}(J^\nn)$ by (2.12)
and (3.16), we obtain $\g'_{\l_1,\l_2}(J^\nn)=0$.\par%
Now suppose $\nn=(0,n_2,n_2+1)$ for $n_2\ge1$. From
$\g'_{\ll[3]}(J^\kk)=0$ for $\kk=(0,0,n_2,n_2+2)$ and
$\kk=(0,0,n_2+1,n_2+1)$, we obtain that
$\g'_{\l_1,\l_2}(J^\nn)=\l_1\g'_{1,\l_2}(J^\nn)$ and that (4.26)
again holds. From this, we obtain that
$\g'_{\l_1,\l_2}(J^\nn)=-\g'_{\l_1,-\l_2}(J^\nn)$. Thus we can
define a reduced $2$-cochain $f'$ such that
$f'_{\l_1}(J^{m_1},J^{m_2})=\g'_{1,\l_1}(J^\nn)$ if
$(m_1,m_2)=(n_2,n_2)$ or $f'_{\l_1}(J^{m_1},J^{m_2})=0$ otherwise.
Then the rest of the proof is as before.\par%
Finally suppose $\nn=(0,n_2,n_3)$ with $n_3\ge n_2+2$ or
$\nn=(n_1,n_2,n_3)$ with $n_1\ge1$. Then the proof is the same as
that of Lemmas 4.4 and 4.5.\qed%
This completes the proof of Theorem 2.5.
\par\vskip-5pt\
\par
\ni{\bf ACKNOWLEDGMENT }%
\par%
This work is supported by a NSF grant 10171064 of China
and two grants ``Excellent Young Teacher Program'' and
``Trans-Century Training Programme Foundation for the Talents''
from Ministry of Education of China.
\par\ \vs{-5pt}\par
\par%
\small \ni\hi2ex\ha1
$^{1\ \,}$C.~Boyallian, V.~G.~Kac and J.~I.~Liberati, ``On the
classification of subalgebras of Cend$_N$ and $gc_N$,''
J.~Algebra {\bf260}, 32--63 (2003).
\par\ni\hi2ex\ha1
$^{2\ \,}$B.~Bakalov, V.~G.~Kac and A.~A.~Voronov, ``Cohomology of
conformal algebras,'' Commun.~Math. Phys. {\bf200}, 561--598
(1999).
\par\ni\hi2ex\ha1
$^{3\ \,}$S.~J.~Cheng and V.~G.~Kac, ``Conformal Modules,'' Asian
J.~Math. {\bf1}, 181--193  (1997). Erratum, Asian J.~Math.
{\bf2}, 153--156 (1998). 
\par\ni\hi2ex\ha1
$^{4\ \,}$S.~J.~Cheng, V.~G.~Kac and M.~Wakimoto, ``Extensions of
conformal modules,'' in {\it Topological field theory, primitive
forms and related topics,} Proceedings of Taniguchi and RIMS
symposia, Progress in Math., Birkh\"auser, 1998.
\par\ni\hi2ex\ha1
$^{5\ \,}$A.~D'Andrea and V.~G.~Kac, ``Structure theory of finite
conformal algebras,'' Selecta Math. {\bf4}, 377--418 (1998).
\par\ni\hi2ex\ha1
$^{6\ \,}$A.~De Sole and V.~G.~Kac, ``Subalgebras of $gc_N$ and
Jacobi polynomials,'' Canad.~Math.~Bull. {\bf 45}, 567--605
(2002).
\par\ni\hi2ex\ha1
$^{7\ \,}$B.~L.~Feigin, ``On the cohomology of the Lie algebra of
vector fields and of the current algebra,'' Selecta Math.~Soviet.
{\bf7}, 49--62 (1988).
\par\ni\hi2ex\ha1
$^{8\ \,}$B.~L.~Feigin and D.~B.~Fuchs, ``Homology of the Lie
algebra of vector fields on the line,'' (Russian) Funkc.~Anal.~i
Pril. {\bf14}, 45--60 (1980).
\par\ni\hi2ex\ha1
$^{9\ \,}$D.~B.~Fuchs, {\it Cohomology of infinite-dimensional Lie
algebras}, Contemporary Soviet Mathematics. Consultants Bureau,
New York, 1986.
\par\ni\hi2ex\ha1
$^{10\,}$D.~Fattori and V.~G.~Kac, ``Classification of finite
simple Lie conformal superalgebras,'' preprint, math.QA/0106002
\par\ni\hi2ex\ha1
$^{11\,}$I.~M.~Gelfand and D.~B.~Fuchs, ``Cohomologies of the Lie
algebra of formal vector fields,'' (Russian) Izv.~Akad.~Nauk SSSR
Ser.~Mat. {\bf34}, 322--337 (1970). 
\par\ni\hi2ex\ha1
$^{12\,}$V.~G.~Kac, {\it Vertex algebras for beginners},
University Lecture Series, 10. American Mathematical Society,
Providence, RI, 1996.
\par\ni\hi2ex\ha1
$^{13\,}$V.~G.~Kac, ``The idea of locality,'' in {\it Physical
applications and mathematical aspects of geometry, groups and
algebras,} H.-D.
Doebner et al, eds., World Sci., Singapore, 16--32 (1997). 
\par\ni\hi2ex\ha1
$^{14\,}$V.~G.~Kac, ``Formal distribution algebras and conformal
algebras,'' a talk at the Brisbane, in {\it Proc.~XIIth
International Congress of Mathematical Physics (ICMP '97)
(Brisbane),} 80--97.
\par\ni\hi2ex\ha1
$^{15\,}$B.~Kostant, ``Lie algebra cohomology and the generalized
Borel-Weil theorem,'' Annals of Math. {\bf74}, 329--387 (1961).
\par\ni\hi2ex\ha1
$^{16\,}$W.~Li, ``2-Cocycles on the algebra of differential
operators,'' J.~Algebra {\bf122}, 64--80 (1989).
\par\ni\hi2ex\ha1
$^{17\,}$J.~L.~Loday, {\it Cyclic homology,} Grundlehren der
Mathematischen Wissenschaften, 301. Springer-Verlag, Berlin, 1992.
\par\ni\hi2ex\ha1
$^{18\,}$Y.~Su, ``2-Cocycles on the Lie algebras of generalized
differential operators,'' Commun.~Alg. {\bf30}, 763--782 (2002).
\par\ni\hi2ex\ha1
$^{19\,}$Y.~Su and K.~Zhao, ``Second cohomology group of
generalized Cartan type W Lie algebras and central extensions,''
Commun.~Alg. {\bf30}, 3285--3309 (2002).
\par\ni\hi2ex\ha1
$^{20\,}$X.~Xu, ``Equivalence of conformal superalgebras to
Hamiltonian superoperators,'' Algebra Colloquium {\bf8}, 63--92
(2001).
\par\ni\hi2ex\ha1
$^{21\,}$X.~Xu, ``Simple conformal algebras generated by Jordan
algebras,'' preprint, math.QA/0008224. 
\par\ni\hi2ex\ha1
$^{22\,}$X.~Xu, ``Simple conformal superalgebras of finite
growth,'' Algebra Colloquium {\bf7}, 205--240 (2000).
\par\ni\hi2ex\ha1
$^{23\,}$X.~Xu, Quadratic Conformal Superalgebras, {\it
J.~Algebra} {\bf231}, 1--38 (2000).
\end{document}